\documentclass[10pt]{article}
\usepackage{amsmath,amsthm,amsfonts,amssymb}
\usepackage{epsfig}
\usepackage{color}
\usepackage{latexsym}
\usepackage{supertabular}
\usepackage{graphicx}
\usepackage{a4wide}
\numberwithin{equation}{section}

\def\whitebox{{\hbox{\hskip 1pt
 \vrule height 6pt depth 1.5pt
 \lower 1.5pt\vbox to 7.5pt{\hrule width
    3.2pt\vfill\hrule width 3.2pt}%
 \vrule height 6pt depth 1.5pt
 \hskip 1pt } }}
\def\qed{\ifhmode\allowbreak\else\nobreak\fi\hfill\quad\nobreak
     \whitebox\medbreak}

\newcommand{\ignore}[1]{}

\theoremstyle{plain}
\newtheorem{theorem}{Theorem}[section]

\newtheorem{example}[theorem]{Example}

\newtheorem{remark}[theorem]{Remark}

\def\qed{{\hfill$\square$}}
\def\proof{{\vspace{-0.3cm}\bf Proof: \,}}

\def\Z{{\mathbb Z}}
\def\Q{{\mathbb Q}}

\def\C{{\mathbb C}}
\def\F{{\mathbb F}}
\def\mod{{\mathrm{mod\,\,}}}

\def\Tr{{\mathrm{Tr}}}
\def\Norm{{\mathrm{Norm}}}

\def\ord{{\mathrm{ord}}}

\title{Evaluation of the weight distribution of a class of cyclic codes based on  index $2$ Gauss sums}
\author{Tao Feng\footnotemark[1] \,\, and\,\,  Koji Momihara\footnotemark[2]}
\date{} 
\begin{document}
\maketitle
\footnotetext[1]{
Department of Mathematics, Zhejiang University, Hangzhou 310027, Zhejiang, China; Email address: 
tfeng@zju.edu.cn}
\footnotetext[2]{
Department of Mathematics, Faculty of Education, Kumamoto University,  
2-40-1 Kurokami, Kumamoto 860-8555, Japan; Email address: 
momihara@educ.kumamoto-u.ac.jp}
\renewcommand{\thefootnote}{\arabic{footnote}}
\begin{abstract}  
The duals of cyclic codes with two zeros have been extensively studied, 
and their weight distributions have recently been evaluated in some cases (\cite{D12,MZ11,WTQYX11,X12}). 
In this note, we determine the weight distribution of a certain new class of such codes by  computations involving index $2$ Gauss sums.
\end{abstract}
\begin{center} 
\end{center}
\section{Introduction}
Let $p$ be a prime, $q$ be a power of $p$, and  $k$ be a positive integer. An {\it $[n,k,d]$-linear code} ${\mathcal C}$ is a $k$-dimensional subspace of $\F_q^n$ with
minimum distance $d$. Each element in ${\mathcal C}$ is called a {\it codeword}. If any cyclic shift of each codeword in ${\mathcal C}$ is again in ${\mathcal C}$, the code is called {\it cyclic}.
To determine all the nonzero weights and their frequencies of a given  code is one of  the main problems in algebraic coding theory.
The {\it weight enumerator} of ${\mathcal C}$ is defined as the polynomial
$
1+\sum_{i=1}^{n}A_i x^i$,
where $A_i$ is the number of codewords of weight $i$ in ${\mathcal C}$. Furthermore,
the sequence $(A_1,A_2,\ldots,A_n)$ is called the {\it weight distribution} of the code
${\mathcal C}$.
Many important families of cyclic codes
have been extensively studied  in the literature, but the weight distributions are generally difficult to compute and there are only a few special families that this has been done.  We assume that the reader is familiar with the basic facts about coding theory, see for instance  \cite{MS06} and \cite{S09}.

Let $\alpha$ be a primitive element of $\F_{q^k}$, and let $h$ and $e$ be positive integers such that $e\,|\,h$ and $h\,|\,q-1$. Put
$g=\alpha^{(q-1)/h}$, $\beta=\alpha^{(q^k-1)/e}$, and $n=h(q^k-1)/(q-1)$. Since the order of $g^{-1}$ and $(\beta g)^{-1}$ are both equal to $n$,
the minimal polynomials $f_1(x)$ and $f_2(x)$ of $g^{-1}$ and $(\beta g)^{-1}$ divide $x^n-1$. Furthermore, it is easy to show that $g^{-q^j} \ne(\beta g)^{-1}$ for any integer $j$,  so we have $f_1(x)f_2(x)\,|\,x^n-1$. By  Delsart's Theorem \cite{D75}, the cyclic code ${\mathcal C}_{(q,k,h,e)}$ with
$f_1(x)f_2(x)$ as its parity-check polynomial can be represented in the following trace form. Let
\[
{\bf c}(a,b)=(\Tr_{q^k/q}(a g^0+b(\beta g)^0),\Tr_{q^k/q}(a g^1+b(\beta g)^1),\ldots,\Tr_{q^k/q}(a g^{n-1}+b(\beta g)^{n-1})),
\]
where $\Tr _{q^k/q}$ is the relative trace from $\F_{q^k}$ to $\F_q$.
Then,
it holds that
\[
{\mathcal C}_{(q,k,h,e)} =\{{\bf c}(a,b)\,|\,(a,b)\in \F_{q^k}^2\}.
\]
The dimension of the code is  $2k$. The recent interest in the weight distribution of this type of codes ${\mathcal C}_{(q,k,h,e)}$ starts with \cite{MZ11}, and is followed by \cite{D12}, \cite{WTQYX11}, \cite{X12}. The objective of this note is to compute the weight distribution of a further class of such codes.

In general, to evaluate the weight distribution of the code ${\mathcal C}_{(q,k,h,e)}$ is quite difficult and most cases  remain unsettled. The weight distribution appears  mostly rather complicated, but still there are cases where not so many nonzero weights are involved and a neat expression is available. Here we list the known cases in which the weight distributions have been explicitly evaluated.
\begin{itemize}
\item[(i)] $e>1$ and $m=1$ \cite{MZ11},
\item[(ii)] $e=2$ and $m=2$ \cite{MZ11},
\item[(iii)] $e=2$ and $m=3$ \cite{D12},
\item[(iv)] $e=2$ and $-1\in \langle p\rangle \,(\mod{m})$ \cite{D12},
\item[(v)] $e=3$ and $m=2$ \cite{WTQYX11},
\item[(vi)] $e=4$ and $m=2$ \cite{X12},
\end{itemize}
where $m=\gcd{(\frac{q^k-1}{q-1},\frac{e}{h}(q-1))}$.
Furthermore, if we set $h=q-1$ and drop the condition $e|h$, then they are related to primitive cyclic codes with two zeros and have been extensively studied in the literature, see for instance \cite{BM10,C04,M04,MR07,S95,YCD06} and the references therein.

The purpose of this note is to compute the weight distribution of ${\mathcal C}_{(q,k,h,e)}$ for the case where $e=2$, $m$ is a prime, the subgroup $\langle p\rangle$ generated by $p\in \Z_{m}^\ast$ has index $2$ in $\Z_{m}^\ast$ and $-1\not\in \langle p\rangle $. Our evaluation is based on the explicit determination of certain index $2$ Gauss sums and the Davenport-Hasse theorem.

\section{Index $2$ Gauss sums }
Let $p$ be a prime, $f$ a positive integer, and $q=p^f$. The canonical additive character $\psi$ of $\F_q$ is defined by
$$\psi\colon\F_q\to \C^{*},\qquad\psi(x)=\zeta_p^{\Tr _{q/p}(x)},$$
where $\zeta_p={\rm exp}(\frac {2\pi i}{p})$, and $\Tr _{q/p}$ is the absolute trace. For each multiplicative character
$\chi$ of $\F_q^\ast$, we define a {\it Gauss sum} over $\F_q$ as follows:
\[
G_q(\chi)=\sum_{x\in \F_q^\ast}\chi(x)\psi(x).
\]

Below are a few basic properties of Gauss sums \cite{LN97}:
\begin{itemize}
\item[(i)] $G_q(\chi)\overline{G_q(\chi)}=q$ if $\chi$ is nontrivial;
\item[(ii)] $G_q(\chi^p)=G_q(\chi)$;
\item[(iii)] $G_q(\chi^{-1})=\chi(-1)\overline{G_q(\chi)}$;
\item[(iv)] $G_q(\chi)=-1$ if $\chi$ is principal.
\end{itemize}

In general, the explicit evaluation of Gauss sums is a very difficult problem. There are only a few cases where the Gauss sums have been evaluated.
The most well known case is the {\it quadratic} case where the order of $\chi$ is two. In this case, it holds that
\begin{equation}\label{eq:quad}
G_{p^f}(\chi)=(-1)^{f-1}\left(\sqrt{p^\ast}\right)^f,\;p^\ast=(-1)^{\frac{p-1}{2}}p.
\end{equation}
The next well-studied case is the so-called {\it semi-primitive } case, where there
exists an integer $j$ such that $p^j\equiv -1\,(\mod{N})$, with $N$ being the order of
 $\chi$. Please refer to \cite{BEW97,BMW82,CK86} for details on the explicit evaluation of Gauss sums in this case.

The next interesting case is the index $2$ case, where the subgroup $\langle p\rangle$ generated by $p\in \Z_{N}^\ast$ has index $2$ in $\Z_{N}^\ast$ and $-1\not\in \langle p\rangle $. In this case,
it is known that $N$ can have at most two odd prime divisors.
Many authors have investigated this case, see e.g., \cite{BM73,L97,M98,MV03,YX10}. In particular, a complete solution to the problem of explicitly evaluating Gauss sums in this case is recently given in \cite{YX10}. We record here the following result which we shall use in the next section.
\begin{theorem}\label{Sec2Thm1}(\cite{YX10},  Case A; Theorem~4.1)
Let $N=p_1^\ell$, where $p_1$  is a prime $\equiv 3\,(\mod{4})$ with $p_1>3$. Assume that $p$ is a prime such that
$\ord_{p_1^\ell}(p)=\phi(p_1^\ell)/2$.
Let $f=\phi(N)/2$, $q=p^f$, and $\chi$ be a multiplicative character of order $N$ of $\F_q^\ast$. Then, for
$0\le s\le \ell-1$, we have
\begin{eqnarray*}
G_q(\chi^{p_1^s})&=&p^{\frac{f-cp_1^s}{2}}
\left(\frac{a+b\sqrt{-p_1}}{2}\right)^{p_1^s},
\end{eqnarray*}
where $c$ is the class number of $\Q(\sqrt{-p_1})$, and $a$ and $b$ are integers
determined by $a,b\not\equiv 0\,(\mod{p})$, $4p^{c}=a^2+p_1b^2$,  and $a \equiv -2p^{\frac{f+c}{2}}\,(\mod{p_1})$.
\end{theorem}

Here, we should remark that index $2$ Gauss sums have been successfully applied to the determination of the weight distribution of certain {\it irreducible cyclic codes} in
\cite{BM73}.  Also, recently they have been used in the construction of new infinite families of combinatorial
configurations, such as strongly regular graphs, skew Hadamard difference sets, and association schemes with nice properties (\cite{FX111,FX112,FX113,FMX11}).

To obtain our main result, we will need the following theorems, the first known as the Davenport-Hasse theorem.
\begin{theorem}\label{thm:lift}(\cite[Theorem 5.14]{LN97})
Let $\chi$ be a nonprincipal multiplicative character on $\F_q^\ast=\F_{p^f}^\ast$ and
let $\chi'$ be the lifted character of $\chi$ to $\F_{q^s}^\ast $, i.e., $\chi'(\alpha):=\chi(\Norm_{\F_{q^s}/\F_q}(\alpha))$ for $\alpha\in \F_{q^s}$.
Then, it holds that
\[
G_{q^{s}}(\chi')=(-1)^{s-1}(G_{q}(\chi))^s.
\]
\end{theorem}
\begin{theorem}(\cite[Theorem 5.30]{LN97})
Let $\psi$ be the canonical additive character of $\F_q$ and
$\chi$ be a multiplicative character of $\F_q$ of order $d\,|\,q-1$.
Then, it holds that
\[
\sum_{x\in \F_q}\psi(ax^d+b)=\psi(b)\sum_{i=1}^{d-1}\chi^{-i}(a)G_q(\chi^i)
\]
for any $a,b\in \F_q$ with $a\not=0$.
\end{theorem}
\section{The weight distribution}
In this section, we shall use the same notations as in the Introduction. Moreover, we fix the settings as follows:
\begin{itemize}
\item[(i)] $e=2$;
\item[(ii)] $m=\gcd{(\frac{q^k-1}{q-1},\frac{e}{h}(q-1))}$ is a prime $\equiv 3 \,(\mod{4})$, which we write as $p_1$;
\item[(iii)] $q=p^f$, $fk$ is divisible by $\frac{p_1-1}{2}$, say $fk=s\frac{p_1-1}{2}$ for some positive integer $s$;
\item[(iv)] $p$ is of  index $2$ modulo $p_1$.
\end{itemize}
Under these assumptions, we will determine the weight distribution of the cyclic code ${\mathcal C}_{(q,k,h,e)}$.

Let $\alpha$ be a fixed primitive element of $\F_{q^k}$, and for each $x\in\F_{q^k}$,. We define $C_i^{(\ell,q^k)}:=\alpha^{i}\langle \alpha^{\ell}\rangle$ for  any $\ell\,|\,q^k-1$, $i\in\Z$. Then,
for any $a,b\in \F_{q^k}$, the Hamming weight of ${\bf c}(a,b)$ is $n-Z(q^k,a,b)$, where
\[
Z(q^k,a,b)=|\{x\in C_0^{((q-1)/h,q^k)}\,|\,\Tr_{q^k/q}(ax+\beta^{\log_\alpha(x)}bx)=0\}|,
\]
where we use
From \cite{D12,MZ11}, we have the following formula on $Z(q^k,a,b)$:
\begin{equation}\label{eq:weight1}
Z(q^k,a,b)=\frac{h(q^k-1)}{q(q-1)}+\frac{h}{eq} m
\sum_{i=0}^{e-1}
\sum_{x\in C_{(q-1)i/h}^{(m,q^k)}}\psi((a+\beta^ib)x),
\end{equation}
where $\psi$ is the canonical additive character of $\F_{q^k}$.
\begin{remark}\label{re1}
By Theorems~\ref{Sec2Thm1} and \ref{thm:lift},
the Gauss sum  $G_{q^k}(\chi)$ with $\chi$ a multiplicative character of order $p_1$ of $\F_{q^k}$ is given as
\begin{eqnarray*}
G_{q^k}(\chi)&&=(-1)^{s-1}\left(G_{p^{(p_1-1)/2}}(\chi')\right)^s\\
  &&=(-1)^{s-1}p^{\frac{(p_1-1-2c)s}{4}}\left(\frac{a+b\sqrt{-p_1}}{2}\right)^{s}\in\Q(\sqrt{-p_1}),
\end{eqnarray*}
where $a,b,$ and $c$ are as defined in Theorem~\ref{Sec2Thm1}, and $\chi'$ is a character of $\F_{p^{(p_1-1)/2}}$ whose lift to $\F_{q^k}^\ast$ is $\chi$.
To ease the notation, we introduce the integers $a_s$, $b_s$ such that
\[
\frac{a_s+b_s\sqrt{-p_1}}{2}:=\left(\frac{a+b\sqrt{-p_1}}{2}\right)^{s}.
\]
\end{remark}
 We comment that we allow $b_s$ to have a sign ambiguity of $\pm 1$. We are now ready to prove our main result.
\begin{theorem}
Let ${\mathcal C}_{(q,k,h,e)}$ be the  $[n,2k]$ cyclic code satisfying the above assumptions (i)--(iv). Each codeword ${\bf c}(a,b)$ has weight $n-Z(q^k,a,b)$, and we associate to it the number
\[
Y(q^k,a,b):=\frac{eq}{h}(Z(q^k,a,b)-\frac{h(q^k-1)}{q(q-1)})+2.
\]
Then the multiset $\{Y(q^k,a,b)\,|\,a,\,b\in\F_{q^k}\}$ has values and corresponding multiplicities as listed  in Table~\ref{Tab1}.
\begin{table}[h]
\caption{
\label{Tab1}
The values of $Y(q^k,a,b)$ and their corresponding multiplicities
}
$$
\begin{array}{|c|c|}
\hline
\mbox{$Y(q^k,a,b)$}&\mbox{frequency}\\
\hline
2q^k&1\\ \hline
(-1)^sp^{\frac{s(p_1-1-2c)}{4}}(a_s-b_sp_1)&\left(\frac{p_1-1}{2}\right)^2\left(\frac{q^k-1}{p_1}\right)^2\\ \hline
(-1)^sp^{\frac{s(p_1-1-2c)}{4}}(a_s+b_sp_1)&\left(\frac{p_1-1}{2}\right)^2\left(\frac{q^k-1}{p_1}\right)^2\\ \hline
(-1)^sp^{\frac{s(p_1-1-2c)}{4}}(1-p_1)a_s&\left(\frac{q^k-1}{p_1}\right)^2\\ \hline
(-1)^sp^{\frac{s(p_1-1-2c)}{4}}\frac{a_s-b_sp_1}{2}+q^k&\frac{(p_1-1)(q^k-1)}{p_1}\\ \hline
(-1)^sp^{\frac{s(p_1-1-2c)}{4}}\frac{(a_s+b_sp_1)}{2}+q^k&\frac{(p_1-1)(q^k-1)}{p_1}\\ \hline
(-1)^sp^{\frac{s(p_1-1-2c)}{4}}\frac{1-p_1}{2}a_s+q^k&\frac{2(q^k-1)}{p_1}\\ \hline
(-1)^s p^{\frac{s(p_1-1-2c)}{4}}a_s&\frac{(p1-1)^2}{2}\left(\frac{q^k-1}{p_1}\right)^2\\ \hline
\frac{(-1)^s}{2} p^{\frac{s(p_1-1-2c)}{4}}(-a_s(-2+p_1)-b_s p_1)&(p_1-1)\left(\frac{q^k-1}{p_1}\right)^2\\ \hline
\frac{(-1)^s}{2} p^{\frac{s(p_1-1-2c)}{4}} (-a_s(-2+p_1)+b_s p_1)&(p_1-1)\left(\frac{q^k-1}{p_1}\right)^2\\
\hline
\end{array}
$$
\end{table}
\end{theorem}
\proof
By the equation (\ref{eq:weight1}) above, it suffices to compute
the sum
\begin{equation}\label{eq:main}
\sum_{i=0,1}\sum_{z\in C_{(q-1)i/h}^{(p_1,q^k)}}\psi((a+(-1)^ib)x).
\end{equation}
Let $E=\{0,1\}$ and $E_{0}^{a,b}=\{i\in E\,|\,a+(-1)^ib=0\}$. If $i\in E_{0}^{a,b}$, then the inner sum is $(q^k-1)/p_1$. Therefore, we have
\begin{eqnarray*}
& &\sum_{i=0,1}\sum_{z\in C_{(q-1)i/h}^{(p_1,q^k)}}\psi((a+(-1)^i b)z)-\frac{q^k-1}{p_1}|E_0^{a,b}|\\
&=&\frac{1}{p_1}\sum_{i\in E\setminus E_0^{a,b}}\sum_{z\in \F_{q^k}^\ast}\psi((a+(-1)^ib) \alpha^{(q-1)i/h} z^{p_1})\\
&=&\frac{1}{p_1}\sum_{i\in E\setminus E_0^{a,b}}\left(\sum_{z\in \F_{q^k}}\psi((a+(-1)^ib) \alpha^{(q-1)i/h} z^{p_1})-1\right)\\
&=&-\frac{e-|E_0^{a,b}|}{p_1}+\frac{1}{p_1}\sum_{i\in E\setminus E_0^{a,b}}\sum_{j=1}^{p_1-1}\chi^{-j}(a+(-1)^ib)G_{q^k}(\chi^j),
\end{eqnarray*}
where $\chi$ is a fixed multiplicative character of order $p_1$. It follows that
\[
Y(q^k,a,b)=q^k|E_0^{a,b}|+\sum_{i\in E\setminus E_0^{a,b}}\sum_{j=1}^{p_1-1}\chi^{-j}(a+(-1)^ib)G_{q^k}(\chi^j).
\]

Now,
by Remark~\ref{re1},
the Gauss sum  $G_{q^k}(\chi)$ is written as
\[
G_{q^k}(\chi)=(-1)^{s-1}p^{\frac{(p_1-1-2c)s}{4}}\left(\frac{a_s+b_s\sqrt{-p_1}}{2}\right).
\]
Since $G_{q^k}(\chi^p)=G_{q^k}(\chi)$ and
$G_{q^k}(\chi^{-1})=\chi(-1)\overline{G_{q^k}(\chi)}=\overline{G_{q^k}(\chi)}$, the second summand in $Y(q^k,a,b)$ is equal to
\begin{eqnarray}\label{eq}
& &\sum_{i\in E\setminus E_0^{a,b}}\left(G_{q^k}(\chi)\sum_{j\in \langle p\rangle}\chi^{-j}(a+(-1)^ib)
+G_{q^k}(\chi^{-1})\sum_{j\in -\langle p\rangle}\chi^{-j}(a+(-1)^ib)\right)\nonumber \nonumber\\
&=&(-1)^{s-1}p^{\frac{(p_1-1-2c)s}{4}}\sum_{i\in E\setminus E_0^{a,b}}2Re\left\{\left(\frac{a_s+b_s\sqrt{-p_1}}{2}\right)\sum_{j\in \langle p\rangle}\psi_{p_1}(-j\ell_{a+(-1)^ib})\right\},\nonumber
\end{eqnarray}
where $\psi_{p_1}$ is the canonical additive character of $\F_{p_1}$ and  $\ell_{a+(-1)^i b}$ is the integer such that
\[
{\ell_{a+(-1)^i b}}\equiv \log_\alpha (a+(-1)^i b)\,(\mod{p_1}).
\]
Now, we compute the sum
$\sum_{j\in \langle p\rangle}\psi_{p_1}(jx)$. If $x\equiv 0\,(\mod{p_1})$, it is clear that
\[
\sum_{j\in \langle p\rangle}\psi_{p_1}(jx)=\frac{p_1-1}{2}.
\]
Let $\eta$ be the quadratic character of $\F_{p_1}^\ast$.
If $x\not\equiv 0\,(\mod{p_1})$, by (\ref{eq:quad}), it holds that
\begin{eqnarray*}
\sum_{j\in \langle p\rangle}\psi_{p_1}(jx)
&=&
\frac{1}{2}\sum_{j\in \F_{p_1}^\ast}(1+\eta(j))\psi_{p_1}(xj)\\
&=&
\frac{-1+\eta(x)G_{p_1}(\eta)}{2}=
\frac{-1+\eta(x)\sqrt{-p_1}}{2}.
\end{eqnarray*}
Then,
the equation (\ref{eq}) is reformed as
\begin{eqnarray}\label{eq:main2}
& &(-1)^{s-1}2p^{\frac{(p_1-1-2c)s}{4}}Re\left\{\left(\frac{a_s+b_s\sqrt{-p_1}}{2}\right)\right.\\
& &\hspace{3cm}\left.\times
\left(N_0\frac{-1-\sqrt{-p_1}}{2}+N_1\frac{-1+\sqrt{-p_1}}{2}+N_2\frac{p_1-1}{2}\right)\right\}\nonumber
\end{eqnarray}
where $N_0$ and $N_1$ are the numbers of nonzero squares and nonsquares modulo $p_1$ in $\{\ell_{a+(-1)^{i}b}\,|\,i\in E\setminus E_0^{a,b}\}$, respectively, and $N_2$ is the number of zeros modulo $p_1$ in $\{\ell_{a+(-1)^{i}b}\,|\,i\in E\setminus E_0^{a,b}\}$.  After simplification, we see that the expression in (\ref{eq}) is equal to $\frac{(-1)^{s}}{2}p^{\frac{(p_1-1-2c)s}{4}}$ times
\[
(N_0+N_1+N_2-p_1N_2)a_s+(N_1-N_0)p_1b_s.
\]

 Since $e=2$, there are ten possibilities for the values of the tuple ($|E_0^{a,b}|,N_0,N_1,N_2$). We shall compute the frequency of each plausible tuple ($|E_0^{a,b}|,N_0,N_1,N_2$), which we denote by $N_{|E_0^{a,b}|,N_0,N_1,N_2}$.  As  a consequence, we obtain the values and multiplicities of $Y(q^k,a,b)$ in Table~\ref{Tab1}.

 It is clear that $N_{2,0,0,0}=1$, $N_{1,1,0,0}=N_{1,0,1,0}=\frac{(p_1-1)(q^k-1)}{p_1}$, $N_{1,0,0,1}=\frac{2(q^k-1)}{p_1}$. For instance, $N_{1,1,0,0}$ is the number of pairs $(a,b)$ such that $a+b=0$, $\log_\alpha(a-b)=\log_\alpha(2a)\pmod{p_1}$ is a nonzero square or $a-b=0$, $\log_\alpha(a+b)=\log_\alpha(2a)\pmod{p_1}$ is a nonzero square, which is easily seen to be $2\cdot\frac{p_1-1}{2}\cdot\frac{q^k-1}{p_1}=\frac{(p_1-1)(q^k-1)}{p_1}$. The other three numbers are obtained similarly.

  Furthermore, by considering the case where at least one of $a+b,\,a-b$ is in $ \bigcup_{i\in \langle p\rangle}C_i^{(p_1,q^k)}$, we have
\begin{eqnarray*}
& &N_{1,1,0,0}+2N_{0,2,0,0}+N_{0,1,1,0}+N_{0,1,0,1}\\
&=&|\{(a,b)\in \F_{q^k}^2\,|\,a+b\in \bigcup_{i\in \langle p\rangle}C_i^{(p_1,q^k)}\}|
+|\{(a,b)\in \F_{q^k}^2\,|\,a-b\in \bigcup_{i\in \langle p\rangle}C_i^{(p_1,q^k)}\}|\\
&=&q^k\frac{(q^k-1)(p_1-1)}{p_1}.
\end{eqnarray*}
Similarly, we have
\begin{eqnarray*}
N_{1,0,1,0}+N_{0,1,1,0}+2N_{0,0,2,0}+N_{0,0,1,1}&=&q^k\frac{(q^k-1)(p_1-1)}{p_1},\\
N_{1,0,0,1}+N_{0,1,0,1}+N_{0,0,1,1}+2N_{0,0,0,2}&=&2q^k\frac{q^k-1}{p_1}.
\end{eqnarray*}

It is therefore enough to compute the values of $N_{0,2,0,0}$, $N_{0,0,2,0}$, and
$N_{0,0,0,2}$ only. Now, we introduce  the following notations: for $b\not=0$
\begin{eqnarray*}
u&:=&ab^{-1}\in \F_{q^k}\setminus\{\pm 1\},\\
t&:=&\log_{\alpha}(u+1)\pmod{p_1}.\\
s&:=&\log_{\alpha}(u-1)\pmod{p_1},\\
x&:=&\log_{\alpha}(b)\pmod{p_1},
\end{eqnarray*}
and
\[
M:=\{u\in \F_{q^k}\setminus \{\pm 1\}\,|\,(u+1)/(u-1)\in C_{0}^{(p_1,q^k)}\}.
\]
Note that
$|M|=(q^k-1)/p_1-1$. Moreover, we will use the well known fact
(\cite[p.~81]{BEW97}) that
\begin{eqnarray*}
& &|(\langle p\rangle +u)\,(\mod{p_1}) \cap \langle p\rangle \,(\mod{p_1})|\\
&=&|(-\langle p\rangle +u) \,(\mod{p_1}) \cap -\langle p\rangle\,(\mod{p_1}) |
=\frac{p_1-3}{4}.
\end{eqnarray*}

Now we are ready to compute the values of $N_{0,2,0,0}$, $N_{0,0,2,0}$, and
$N_{0,0,0,2}$, from which all the remaining numbers $N_{0,1,1,0}$, $N_{0,1,0,1}$, $N_{0,0,1,1}$ will follow.

(1) {\bf $N_{0,0,0,2}$: } Recall that
\[
N_{0,0,0,2}=|\{(a,b)\in \F_{q^k}^2\,|\,a+b,a-b\in C_0^{(p_1,q^k)}\}|.
\]
There are $(q^k-1)/p_1$ such pairs with $b=0$.  Assume $b\ne 0$. Then $a+b,a-b\in C_0^{(p_1,q^k)}$ amounts to $t+x=0,s+x=0$ in $Z_{p_1}$, so $t=s$, i.e., $\frac{u+1}{u-1}\in C_0^{(p_1,q^k)}$, which amounts to saying that $u\in M$. For each such $u$, there is a unique $x\in \Z_{p_1}$, so $(q^k-1)/p_1$ of $b\in C_x^{(p_1,q^k)}$. Now we have
\[
N_{0,0,0,2}=\frac{q^k-1}{p_1}+\frac{q^k-1}{p_1}\cdot |M|=\left(\frac{q^k-1}{p_1}\right)^2.
\]

(2){ \bf $N_{0,2,0,0}$: } Recall that
\[
N_{0,2,0,0}=|\{(a,b)\in \F_{q^k}^2|a+b,a-b\in \bigcup_{i\in \langle p\rangle}C_i^{(p_1,q^k)}\}|.
\]
There are $\frac{p_1-1}{2}\cdot \frac{q^k-1}{p_1}$ such pairs with $b=0$. Assume $b\ne 0$. Then $a+b,a-b\in \bigcup_{i\in \langle p\rangle}C_i^{(p_1,q^k)}$ amounts to $t+x,s+x\in \langle p\rangle$, i.e., $x\in (\langle p\rangle-t+s)\cap \langle p\rangle$. There are two cases:
(1) $s=t$ and there are $\frac{p_1-1}{2}$ $x$'s, (2) $s\ne t$ and there are $\frac{p_1-3}{4}$ $x$'s. Note that $s=t$ if and only if $u\in M$. A similar argument as in the determination of $N_{0,2,0,0}$ gives that
\begin{eqnarray*}
& &N_{0,2,0,0}\\
&=&\left(\frac{p_1-1}{2}\right)\left(\frac{q^k-1}{p_1}\right)
+|M|\cdot\left(\frac{p_1-1}{2}\right)\left(\frac{q^k-1}{p_1}\right)
+(p_1-1)\left(\frac{p_1-3}{4}\right)\left(\frac{q^k-1}{p_1}\right)^2\\
&=&\left(\frac{p_1-1}{2}\right)^2\left(\frac{q^k-1}{p_1}\right)^2.
\end{eqnarray*}

(3) {\bf $N_{0,0,2,0}$: }Proceed exactly the same way as above and we obtain $N_{0,0,2,0}=\left(\frac{p_1-1}{2}\right)^2\left(\frac{q^k-1}{p_1}\right)^2$.

To sum up, we get the result listed in Table~\ref{Tab1}.
\qed
\begin{example}
Consider the case where
\[
(p,f,k,e,h,m)=(3,5,55,2,2,11).
\]
In this case, $p=3$ is index $2$ modulo $11$ and the class number of
$\Q(\sqrt{-11})$ is $1$.
The Gauss sum $G_{3^{5}}(\chi)$ with $\chi$ a character  of order $11$ of $\F_{3^5}$ is given as $
(1\pm \sqrt{-11})/2$,
where the sign ambiguity $\pm$ will not matter. Then, the Gauss sum  $G_{3^{5\cdot 11}}(\chi')$ for the lifted character $\chi'$  of $\chi$ is given as
\[
\left(\frac{1\pm \sqrt{-11}}{2}\right)^{11}=\frac{67\pm 253\sqrt{-11}}{2},
\]
i.e., $a_{11}=67$, $b_{11}=\pm 253$. By Table~\ref{Tab1}, the code ${\mathcal C}_{q,k,h,e}$ is a $[2(3^{55} - 1)/(3^5 - 1),22, 3^{17}(-1358+3^{33})]$-linear code over $\F_{3^5}$ with the following weight distribution:
\begin{eqnarray*}
& &1+25\ell^2 x^{2A(B-1358)}+25\ell^2 x^{2A(B+1425)}+
\ell^2x^{2A(B-335)}
+10\ell x^{A(B-1358)}+10\ell x^{A(B+1425)}\\
& &+2\ell x^{A(B-335)}+50\ell^2 x^{A(2B+67)}+10\ell^2 x^{A(2B-1693)}+10\ell^2 x^{2A(B+545)},
\end{eqnarray*}
where $A=3^{17},B=3^{33},\ell=(3^{55}-1)/11$.
\end{example}

\section{Conclusion}
In this short note, we explicitly determine the weight distribution of a class of cyclic codes ${\mathcal C}_{(q,k,h,e)}$ under certain index $2$ condition as specified at the beginning of Section 3 when $e=2$.

Under the  assumptions (ii)--(iv), if we allow  $e>1$ to be arbitrary, then there will be $\binom{e+3}{3}$ possible weights, and at least $\binom{e+2}{2}$ of them have roughly the same (nonzero) count when $q^k$ is large compared to $p_1^e$, according to the estimate by Xiong \cite{X12}. For instance, theoretically we should be able to determine the weight enumerator under the  assumptions (ii)--(iv) when $e=3$ using the same technique here by more involved computations, but in general there will be $20$ weights. Therefore, it will be of interest to determine the cases where there are only few nonzero weights, say, less than ten. We leave this for future work.

If we have a multiplicative character $\chi$ of prime order $p_1$ over a finite field $\F_{q_1}$, and $G_{q_1}(\chi)$ is in the quadratic subfield of  $\Q(\zeta_{p_1})$, then our method also apply and yield similar results. In our construction, the index 2 condition is used to guarantee this point.  So it will be interesting to determine all such Gauss sums. We leave this for future work.

 As we see in the application of applying Gauss sums to the construction of combinatorial objects, we first succeed in the index 2 case, and then extend  to the index 4 case and then even more complicated settings. We wonder this will be the case in the application discussed in this note. We leave this for future work.

\section*{Acknowledgment}

\end{document}